\documentclass[12pt]{amsart}
\newcommand{\Prod}{\Pi}
\newcommand{\A}{\mathbb A}

\newcommand{\C}{\mathbb C}

\newcommand{\Aut}{\mathop{Aut}}
\newcommand{\Mor}{\mathop{Mor}}
\newcommand{\dd}[1]{\mathop{\frac{\partial}{\partial #1}}}
\newcommand{\Spec}{\mathop{Spec}}
\newcommand{\Sing}{\mathop{Sing}}

\newcommand{\MSC}[1]{{\renewcommand{\thefootnote}{}%
\footnote{MSC Subject Classification: #1}}}
\newtheorem{theorem}{Theorem}
\newtheorem{proposition}{Proposition}
\newtheorem{lemma}{Lemma}

\newtheorem{corollary}{Corollary}

\newtheorem{definition}{Definition}

\newtheorem*{assumptions}{Key Assumptions}
\newtheorem*{example}{Example}
\newtheorem*{Falsity}{Falsity}
\begin{document}
\title{Invariant Rings and Quasiaffine Quotients.}
\author{J\"org Winkelmann}
\MSC{Primary: 13A50; Secondary: 14R20, 14L30}
\begin{abstract}
We study Hilbert's fourteenth problem from a geometric point of view.
Nagata's celebrated counterexample demonstrates that for an arbitrary
group action on a variety the ring of invariant functions need not
be isomorphic to the ring of functions of an affine variety.
In this paper we will show that nevertheless it is always isomorphic
to the ring of functions on a quasi-affine variety.
\end{abstract}
\maketitle
\section{Introduction}
The fourteenth of Hilbert's famous problems (\cite{H}) is the following.

\medskip
\begin{em}
Let $K/L$ and $L/k$ be field extensions, and $A\subset K$ be
a finitely generated $k$-algebra.

Does this imply that $A\cap L$ is also a finitely generated $k$-algebra?
\end{em}

\medskip
This problem was motivated by the following special case:

\medskip
\begin{em}
Let $k$ be a field and $G\subset GL(n,k)$ a subgroup. Is the ring
of invariants $k[x_1,\ldots,x_n]^G$ a finitely generated $k$-algebra?
\end{em}

\medskip
(This is a special case: Take $K=k(x_1,\ldots,x_n)$
and $L=K^G$.)

For reductive groups this is indeed the case. This was already shown
by Hilbert. However, for non-reductive
groups there is the celebrated counterexample of Nagata
\cite{N1}).
Popov deduced from Nagata's example that for every non-reductive
algebraic group $G$ there exists an affine $G$ variety such that
the ring of invariants is not finitely generated \cite{P}.
In 1990, a new counterexample was found by Roberts \cite{Rob}.
Later, further counterexamples were obtained by
Deveney and Finston \cite{DF} and by
A'Campo-Neuen \cite{ACN}. Recently, Daigle and Freudenburg
constructed examples in dimension $6$ and $5$ (\cite{DgF},\cite{F}).

Reformulated in a more geometric fashion,
Hilbert's 14th problem ask whether the ring of invariant functions
is necessarily isomorphic to the ring of regular functions on some
affine variety.

From this point of view it is maybe not too
surprising that the answer is negative in general. 
Quotients of affine varieties by actions of
(non-reductive) algebraic groups are often quasi-affine
without being affine, and for arbitrary quasi-affine varieties the
ring of regular functions is not necessarily finitely generated
(see e.g.\ \cite{Ne},\cite{Oe},\cite{R}).
Thus even if the ring of invariants is not isomorphic to
the ring of regular functions on an affine variety it nevertheless may
be isomorphic to the ring of regular functions on a quasi-affine
variety.
The purpose of this note is to demonstrate that this is indeed 
always the
case. Actually we show  that the a $k$-algebra occurs as the
ring of invariants for some affine $G$-variety if and only if
it is isomorphic to the algebra of regular functions
on some quasi affine variety.

\begin{theorem}
Let $k$ be a field and $R$ an integrally closed $k$-algebra.

Then the following properties are equivalent:
\begin{enumerate}
\item
There exists an irreducible, reduced
 $k$-variety $V$ and a subgroup $G\subset Aut_k(V)$
such that $R\simeq k[V]^G$.
\item
There exists a quasi-affine irreducible, reduced
$k$-variety $V$ such that
$R\simeq k[V]$.
\item
There exists an affine irreducible, reduced
$k$-variety $V$ and a regular action
of $G_a=(k,+)$ on $V$ defined over $k$ such that
$R\simeq k[V]^{G_a}$.
\end{enumerate}

If $char(k)=0$, these properties are futhermore equivalent
to the following:
\begin{enumerate}
\setcounter{enumi}{3}
\item
There exists a finitely generated, integrally closed
 $k$-algebra $A$ and a
locally nilpotent derivation $D$ on $A$ such that $R\simeq\ker D$.
\end{enumerate}
\end{theorem}

This result is based on the following more general theorem.
\begin{theorem}
Let $k$ be a field, $V$ an irreducible, reduced, normal $k$-variety, and 
$L$ a subfield of the function field $k(V)$,
containing $k$. Let $R=k[V]\cap L$.

Then there exists a finitely generated $k$-subalgebra $R_0$ of $R$
such that
\begin{enumerate}
\item
The quotient fields of $R$ and $R_0$ coincide.
\item
For every prime ideal $p$ of height one in $R$ the prime ideal $p\cap
R_0$ of $R_0$ also has height one.
\item
There is an open $k$-subvariety $\Omega\subset \Spec(R_0)$ such that
$R=k[\Omega]$ (as subsets of $Q(R)$).
\end{enumerate}
\end{theorem}

These results can be used to construct some
``quasi-affine'' quotient for a group action on an algebraic
variety.

\begin{theorem}
Let $k$ be a field, $V$ an irreducible, reduced,
normal $k$-variety and $G\subset Aut(V)$.

Then there exists a quasi-affine $k$-variety $Z$ and a rational map
$\pi:V\to Z$ such that
\begin{enumerate}
\item
The rational map $\pi$ induces an inclusion $\pi^*:k[Z]\subset k[V]$.
\item
The image of the pull-back $\pi^*(k[Z])$ coincides with the ring of
invariant functions $k[V]^G$.
\item
For every affine $k$-variety $W$ and every $G$-invariant morphism
$f:V\to W$ there exists a morphism $F:Z\to W$ such that $F\circ\pi$ is
a morphism and $f=F\circ\pi$.
\end{enumerate}
\end{theorem}

We may also translate our results in the language of category theory
(also know asa ``general nonsense'') and deduce the following.
\begin{theorem}
For a field $k$  let $\mathcal V_k$ denote the category whose objects are
irreducible reduced normal $k$-varieties and whose morphisms are those
dominant rational maps for which the pull-back of every regular
function is again regular. Let $\mathcal Q_k$ denote the full sub-category 
whose objects consist of all \begin{em} quasiaffine
\end{em} such
varieties.

Then for every object $V\in\mathcal V_k$ and every subgroup 
$G\subset \Aut_{\mathcal V_k} (V)$ the
functor $\Mor_{\mathcal V_k}(V,\cdot)^G$ is representable in the category $\mathcal Q_k$.
\end{theorem}
\section{Preparations}
Following ideas of Nagata \cite{N1}
we employ the notions of \begin{em} Krull rings
\end{em} and \begin{em} ideal transforms\end{em} as algebraic
tools for our constructions.
We give proofs for some
basic facts although they are well known. This is for the benefit of being
self-contained
and because the proofs are so short.
\subsection{Krull rings}
\begin{definition}
An integral domain $R$ is called a \begin{em} Krull ring \end{em} if
there is a family $F$ of discrete valuations on the quotient field $K$
of $R$ such that $R=\{f\in K: v(f)\le 0\ \forall v\in F\}$ and
$\{v\in F:v(f)\ne 0\}$ is finite for every $f\in K$.
\end{definition}
For basic facts on Krull rings, see \cite{B}\cite{N1}.
Noetherian integral domains integrally closed in their quotient
fields
are Krull rings. 
Intersections of Krull rings inside a fixed field are again Krull rings.
For any Krull ring the family of valuations $F$ necessarily contains 
valuations associated to all prime ideals of height one, on the other
hand for a Krull ring $F$ can be choosen as the set of all discrete
valuations associated to prime ideals of height one.

Next, we recall that  for every $k$-variety $V$ the ring
of functions $k[V]$ is a Krull ring and for every Krull ring $R$ and
every group $G$ acting on $R$ by ring automorphisms $R^G$ is again a
Krull ring.

\begin{lemma}\label{lem1}
Let $k$ be a field, $V$ an irreducible, reduced and normal $k$-variety.
Then $k[V]$ is a Krull ring.
\end{lemma}
\begin{proof}
Let $(U_i)_{i\in I}$ be a collecting of affine $k$-varieties covering $V$.
For every $i\in I$ the ring $k[U_i]$ is noetherian and integrally
closed 
and therefore a Krull ring.
Now $k[V]$ equals the intersection of all the $k[U_i]$
considered as subrings of the function field $k(V)$. Hence
$k[V]$ is a Krull ring.
\end{proof}
\begin{lemma}\label{lem2}
Let $R$ be a Krull ring and $G\subset Aut(R)$.

Then $R^G$ is a Krull ring.
\end{lemma}
\begin{proof}
This is immediate, because $R^G$ is the intersection of two Krull
rings, namely $R$ and $Q(R)^G$ where $Q(R)$ denotes the quotient field
of $R$.
\end{proof}
\subsection{The Ideal transform}
\begin{definition}
Let $R$ be an integral domain, $I$ an ideal.
Then the \begin{em} $I$-transform \end{em} of $R$ is defined as
$S(I,R)=\{x\in K:\exists n: x(I^n)\subset R\}$ where $K$ denotes the
quotient field of $R$.
\end{definition}
For Krull rings the ideal transform has particular nice descriptions
in geometric as well as in algebraic form.
First we explain the geometric description.
\begin{lemma}\label{lem3}
Let $R$ be a Krull ring, $I$ an ideal.

Then $S(I,R)$ equals the ring of regular functions on
$\Spec R\setminus V(I)$ (both considered as subsets of the quotient field of
R).
\end{lemma}
\begin{proof}
Let $K$ denote the quotient field of $R$. Let $f\in S(I,R)$ and
$p\in X=\Spec R\setminus V(I)$.
Since $p\notin V(I)$, the prime ideal $p$ of $R$ does not contain $I$.
Hence there is an element $h\in I\setminus p$. By definition $fh^n\in
R$ for some natural number $n$. Since $h\notin p$, this implies $f\in
R_p$.

We will now show the converse. Let $f$ be a regular function on
$X$, i.e., $f\in R_p\forall p\not\supset I$, 
and let $\Lambda=\{v\in F:v(f)>0\}$. For $v\in\Lambda$ consider
the associate prime ideal 
$p_v=\{g\in R:v(g)<0\}$. 
Now $v(f)>0$ implies that
$f\notin R_{p_v}$, hence $I\subset p_v$. Thus $v(x)<0$ for all $x\in
I$ and $v\in\Lambda$. Since $\Lambda$ is finite, it follows that
$fI^n\subset R$ for $n$ sufficiently large, i.e.~$f\in S(I,R)$.
\end{proof}

Next we come to the algebraic description of the ideal transform.
\begin{lemma}\label{lem4}
Let $R$ be a Krull ring, $K$ its quotient field, $F$ the set of all discrete
valuations corresponding to prime ideals of height one in $R$
and $I$ an ideal in $R$. Then $S(I,R)=
\{x\in K:v(x)\le 0\ \forall v \in F'\}$ with $F'=\{v\in F:I\not\subset
p_v\}$.
\end{lemma}
\begin{proof}
Let $x\in S(I,R)$ and $v\in F'$. Then there is an element $f\in I$
such that $v(f)=0$. Now $xI^n\subset R$ implies $v(x)\le 0$.
Conversely, assume $x$ is an element in $K$ such that $v(x)\le 0$ for
all $v\in F'$. Let $\Lambda=\{v\in F:v(x)>0$. Then $\Lambda$ is
finite and $I\subset p_v$ for all $v\in\Lambda$. This implies
$xI^n\subset R$ for $n$ sufficiently large.
\end{proof}
We will need the fact that for prime ideals of height one the ideal
transform is non-trivial.
\begin{lemma}\label{lem5}
Let $R$ be a Krull ring, $p$ a prime ideal of height one.
Then $S(p,R)\ne R$.
\end{lemma}
\begin{proof}
Let $K$ be the quotient field of $R$, $R=\{x\in K:v(x)\le 0\ \forall
v\in F\}$, $p=\{x\in R:v_0(x)<0\}$.
Choose $x\in P$ and define $\Lambda=\{v\in F: v(1/x)>0\}\setminus\{v_0\}$. 
Now $\Lambda$ is finite, and $p_v$ being a prime ideal of height one
for every $v$ implies
 that $p_v$ is not contained in $p$ for any $v\in\Lambda$.
Hence there is an element $y\in R$ with $y\not\in p$ but $v(y)<0$ for
all $v\in\Lambda$. Then $y^nx\in S(p,R)\setminus R$ for $n$
sufficiently large.
\end{proof}

\subsection{}%
In the situation in which we are interested prime ideals of height one
are determined by their zero sets.
\begin{lemma}\label{lem6}
Let $k$ be a field, $V$ an integral, normal $k$-variety, $R$ a $k$-sub algebra of
$k[V]$ such that $R$ is a Krull ring and $R=Q(R)\cap k[V]$.

Let $p$ be a prime ideal of height one in $R$.
Then $p=I(Z(p))\cap R$, where $I(Z(p))$ denotes the set of all $f\in
k[V]$ vanishing on the zero set of $p$.

Moreover, $Z(p)$ contains an irreducible component $Z_0$ such that $I(Z_0)=p$.
\end{lemma}
\begin{proof}
Let $v$ denote the discrete valuation corresponding to $p$.
By $(*)$ there is an element $x\in Q(R)\setminus R$ with $v(x)>0$
and a number $N$
such that $xp^N\subset R$. Thus $x$ defines a rational function on $V$
whose poles are contained in $Z(p)$. 
Now assume $g\in I(Z(p))$ but $g\not\in p$. This would imply that
$v(g)=0$ and $g$ vanishes on $Z(p)$. Since the poles of $x$ are
contained in $Z(p)$, it follows that $g^mx\in k[V]$ for $m$ large
enough.
On the other hand $v(g^mx)=v(x)>0$ implies $g^mx\not\in R$. Thus the
assumption $I(Z(p))\ne p$ yields a contradiction to the requirement
that $R=Q(R)\cap k[V]$.

For the final statement, let $(Z_i)_{i\in I}$ denote the irreducible
components of $Z(p)$. Then $I(Z(p))=\cap_i I(Z_i)$. Since both
$I(Z(p))$ and all of the $I(Z_i)$ are prime ideals, it follows that
there exists an $i\in I$ such that $I(Z_i)=I(Z(p))=p$.
\end{proof}

\begin{lemma}\label{lem7}
Let $k$ be a field and $V$ a $k$-variety.
For every $k$-algebra $R\subset k[V]$
let $\sim_R$ denote the equivalence relation
on the geometric ($\bar k$-rational) points of $V$ given
by $x\sim_R y$ iff $f(x)=f(y)$ for all $f\in R$.

Then for every $k$-sub algebra $R\subset k[V]$
there exists a finitely generated $k$-algebra $R_0\subset R$
such that $\sim_R$ coincides with $\sim_{R_0}$

Furthermore $R_0$ can be choosen as a Krull ring, if $R$ is a
Krull ring.
\end{lemma}
\begin{proof}
The equivalence relation $\sim_R$ defines a subset
$E_R\subset V\times V$ via 
\begin{equation}
E_R=\{(x,y)\in V\times V:x\sim_R y\}
=\{(x,y):f(x)=f(y)\,\forall f \in R\}\tag{$*$}
\end{equation}

Then $E_R$ is the  $k$-sub variety defined by the radical
of the ideal of $k[V\times V]$ generated by $\pi_1^*f-\pi_2^*f$ with $f$
running through $R$.
If $R'\subset R$ are two $k$-sub algebras of $k[V]$,
then $E_{R}\subset E_{R'}$.
Since $V\times V$ is noetherien it follows that for any $k$-sub
algebra $R\subset k[V]$ there exists a finitely generated
$k$-sub algebra $R_0\subset R$ such that $E_R=E_{R_0}$.

Finally, if $R$ is a Krull ring, then $R$ is integrally closed
in its quotient field. Hence the integral closure of $R_0$
is again contained in $R$. Furthermore the integral closure
of $R_0$ is again finitely generated as a $k$-algebra. 
Thus we may assume that $R_0$ is integrally closed.
But integrally closed finitely generated $k$-algebras are Krull rings.
\end{proof}

Let us now fix some key assumptions.
\begin{assumptions}
In the sequel, $k$ is a field, $V$ an irreducible, reduced and normal
$k$-variety, $R\subset k[V]$ a $k$-sub algebra such that $Q(R)
\cap k[V]=R$,
$R_0$ is a finitely generated $k$-algebra such that $R_0\subset R$,
$Q(R_0)=Q(R)$, $R_0$ is integrally closed in $Q(R_0)=Q(R)$ and
$E_R=E_{R_0}$ where $E_R$ is defined as in $(*)$ above.
\end{assumptions}

We will prove the statements of theorem 2 hold for 
such a choice of $R_0$.

\begin{lemma}\label{lem8}
Under the key assumptions
 $height(p\cap R_0)=1$ for every prime ideal $p$ of height one in $R$.
\end{lemma}
\begin{proof}
Let $W=\Spec(R_0)$, this is an affine $k$-variety. The prime ideals of
height one in $R_0$ correspond to the irreducible hypersurfaces in
$W$. Let $\tau:V\to W$ denote the morphism induced by $R_0\subset
k[V]$.
Let $p$ be a prime ideal of height one in $R$. We have seen above that
there is an irreducible subvariety $Z\subset V$ such that $I(Z)\cap
R=p$.

From $E_R=E_{R_0}$ we infer that there exists an irreducible
subvariety $Y\subset Z$ such that
$\tau|_Y:Y\to W$ is generically quasi-finite and $Z$ is the smallest
$E_R$-saturated subvariety containing $Y$.
Now let $X$ be an irreducible subvariety of $V$ such that
$\dim(X)=\dim(Y)+1$ and $X\not\subset Z(p)$.
Let $I=I(X)\cap R$. Then $p\ne I$. On the other hand $I\subset I(Y)$
and $I(Y)\cap R=I(Z)\cap R=p$. Thus $height(p)=1$ implies $I=\{0\}$.
It follows that $\tau(X)$ must be Zariski-dense in $W=\Spec(R_0)$.
Since $\dim(X)=\dim(Y)+1$ it now follows from
$\tau|_Y$ being generically quasi-finite that $\overline{\tau(Y)}$
either is a hypersurface or equals $W$. The latter is excluded since
$p\ne\{0\}$ and $Y\subset Z(p)$.
Thus $\overline{\tau(Y)}=\overline{\tau(Z)}$ has to be a hypersurface
implying that $height(p\cap R_0)=1$.
\end{proof}
\begin{lemma}\label{lem9}
Under the key assumptions for any
two distinct prime ideals of height one $p_1,p_2\subset R$ we have
$p_1\cap R_0\ne p_2\cap R_0$.
\end{lemma}
\begin{proof}
This is an immediate consequence of $E_R=E_{R_0}$ and $I(Z(p_i))\cap R=p_i$.
\end{proof}
\begin{corollary}
Let $F$ resp. $F_0$ denote the set of discrete valuations of $K=Q(R)$
corresponding to prime ideals of height one in $R$ resp. $R_0$.
Then $F\subset F_0$.
\end{corollary}
\begin{lemma}\label{lem10}
The set $F_0\setminus F$ is finite.
\end{lemma}
\begin{proof}
Let $v\in F_0\setminus F$, $H_0\subset W=\Spec(R_0)$ the corresponding
hypersurface and $H=\tau^{-1}(H_0)$. We claim that $\tau(H)$ is not
Zariski-dense in $H_0$. Indeed, if it were Zariski-dense, there would
exist an irreducible subvariety $H'$ with $\overline{\tau(H')}=H$.
But this implies $I(H')\cap R_0=I(H_0)$. Now let $p$ be a prime ideal
of height one contained in $I(H')$ (such an ideal exists, because
$I(H')$ is a prime ideal and $R$ is a Krull ring). Then $p\cap R_0$ is
a prime ideal of height one. Since $I(H_0)$ is of height one, it
follows that $I(H_0)=p\cap R_0$ contrary to our assumption $v\not\in
F$.
Thus $H$ is a hypersurface in $W$ with $\tau(\tau^{-1}(H))$ not being
dense in $H$. Since $\tau:V\to W$ is dominant, there are only finitely
many such hypersurfaces in $W$.
\end{proof}
\section{Proofs of the theorems}
\subsection{Proof of theorem 2.}
\begin{proof}
Let $k$, $V$, $L$ and $R$ be as in the theorem.
Note that $L$ is a finitely generated field extension of $k$,
because $k(V)/k$ is finitely generated and $L\subset k(V)$.

By lemma \ref{lem7} there is a finitely generated $k$-algebra $R_1$
with $R_1\subset R$, $E_R=E_{R_1}$ and $R_1$ being a Krull ring.
We may adjoin finitely many further elements of $R$ to $R_1$ 
and thereby assume
that the quotient fields of $R$ and $R_1$ coincide.
Then we choose $R_0$ as the integral closure of $R_1$ in $L$.
Since $R$ is integrally closed, we have $R_0\subset R$.
Furthermore $E_R=E_{R_0}$ and $R_1$ is again finitely generated,
because it is the integral closure of a finitely generated $k$-algebra.

Thus $R_0$ fulfills the ``Key Assumptions''. Statement $(1)$ of the 
theorem is clear, and statement $(2)$ follows from lemma \ref{lem8}.

Next we define $F$ and $F_0$ as in the corollary above.
By lemma \ref{lem10} the difference set $F_0\setminus F$ is finite.
Hence we may define an ideal of $R_0$ by
\[
I=\Prod_{v\in\{F_0\setminus F\}}\, p_v
\]
with $p_v=\{x\in R_0: v(x)<0\}$.

Since $p_\mu$ is of height one for every $\mu\in F_0$, it is clear that
$I\not\subset p_v$ for $v\in F$. 
Therefore $R=S(I,R_0)$ by lemma \ref{lem4}.

Finally, statement $(3)$ of the theorem follows with the aid of
lemma \ref{lem3}.
\end{proof}

Before starting the proof of theorem 3 we need the subsequent lemma.
\begin{lemma}\label{lem11}
Let $k$ be a field, $V$ a $k$-variety and $W$ a quasi-affine
$k$-variety.

Then there is a one-to-one correspondence between $k$-algebra
homomorphisms $\phi:k[W]\to k[V]$ and rational maps $f:V\to W$ with
$f^*k[W]\subset k[V]$.
\end{lemma}
\begin{proof}
Given a $k$-algebra homomorphism  $\phi:k[W]\to k[V]$, choose a
finitely generated $k$-sub algebra $A\subset k[W]$ such that the
quotient fields of $k[W]$ and $A$ coincide.
The restriction of $\phi$ yields a $k$-morphism $F$ from $V$ to
$Z=\Spec(A)$ and the inclusion $A\subset k[W]$ yields a birational
morphism $\tau:W\to Z$. Now $\tau^{-1}\circ F$ is the desired rational
map.
\end{proof}
\subsection{Proof of theorem 3.}
\begin{proof}
We apply theorem~2 with $L=k(V)^G$ and set $Z=\Omega$.
There is an inclusion $k[Z]=k[\Omega]=R\subset k[V]$.
By the preceding lemma this induces a
rational map $\pi:V\to Z$ with $\pi^*k[Z]\subset k[V]$.
Since $R=k[V]\cap L=k[V]^G$, it follows that
$\pi^*k[Z]=k[V]^G$.

Finally note that for every affine $k$-variety $W$ and every
$G$-invariant morphism $f:V\to W$ we obtain an inclusion 
$f^*k[W]\subset k[V]^G\simeq k[Z]$ which implies that there
exists a morphism $F:Z\to W$ such that $f=F\circ\pi$.
\end{proof}

\subsection{Proof of theorem 1.}
\begin{proof}
The implication $(1)\implies(2)$ follows from
theorem~3. $(3)\implies (1)$ is trivial and $(2)\implies(3)$
follows from the proposition below.

Finally, $(3)\iff(4)$ for the case of characteristic zero is implied
by the well-known correspondence between locally nilpotent derivations
and $G_a$-actions on affine varieties in characteristic zero.
\end{proof}

\begin{proposition}
Let $k$ be a field and let $V$ be a normal quasi-affine $k$-variety.

Then there exists a normal affine $k$-variety $W$ and a regular action
of the additive group $G=G_a$ defined over $k$ on $W$ such that
$k[V]\simeq k[W]^G$.
\end{proposition}
\begin{proof}
Let $V\hookrightarrow Y$ be an open embedding in a normal affine
$k$-variety $Y$, and let $S=Y\setminus V$. Let $D$ denote the union
of codimension $1$-components of $S$ and choose a regular function
$f_1$ on $Y$ such that $D$ is contained in the zero set of $f_1$.
Then choose a regular function $f_2$ on $Y$ such that $f_2$
vanishes on $D$, but does not vanish on any irreducible component
of $Z(f_1)\setminus D$. If $char(k)=p>0$, we replace $f_i$ by
$f_i^{p^N}$ for a sufficiently large $N$. In this way we may
assume that both functions $f_i$ are defined over a 
finite Galois extension $k'/k$ with Galois group $\Gamma$.
Now we may replace $f_i$ by $\Prod_{\sigma\in\Gamma}{}^{\sigma}f_i$.
Therefore we may assume that both $f_i$ are defined over $k$.
We obtain a $k$-morphism $f:Y\to \A^2$. By construction $D$ is the
union of codimension $1$-components of $E=f^{-1}\{(0,0)\}$.
Since regular functions extend through subvarieties of codimension
at least 2 on normal varieties, it follows that
\[
k[V]\simeq k[Y\setminus D]\simeq k[\Omega]
\]
for $\Omega=Y\setminus E$.

Next we consider 
\[
S=\left\{
\begin{pmatrix}
a & b \\
c & d 
\end{pmatrix}
: ad-bc=1
\right\}
\]
and the natural projection $\pi:S\to \A^2$ given by
$\pi:(a,b,c,d)\mapsto (a,b)$.
This realizes $\A^2\setminus\{(0,0)\}$ as the quotient
of $S$ by the $G_a$-action given by
\[
t:
\begin{pmatrix}
a & b \\
c & d 
\end{pmatrix}
\mapsto
\begin{pmatrix}
a & b+ta \\
c & d+tc 
\end{pmatrix}
\]
Now the fiber product $W=Y\times_{\A^2} S$ is an affine variety
and $W\simeq \Omega\times_{\A^2\setminus\{(0,0)\}}S$, because
the image of $S$ in $\A^2$ is contained in $\A^2\setminus\{(0,0)\}$
and $\Omega=f^{-1}(\A^2\setminus\{(0,0)\})$.
The $G_a$-action on $S$ induces a $G_a$-action on the fibered
product $W$ and evidently $k[\Omega]\simeq k[W]^{G_a}$.
\end{proof}

\subsection{Proof of theorem~4}
\begin{proof}
Let $X\in\mathcal V_k$ and $G\subset \Aut_{\mathcal V_k}(X)$.
Each element in $G$ is a birational self-map of $X$ such that
the induced field automorphism of $k(X)$ stabilizes $k[X]$.
In particular $G\hookrightarrow \Aut(k[X])$.
Due to theorem~2 there exists an
object $\Omega\in\mathcal Q_k$ and a 
$\mathcal V_k$-morphism $\pi:X\to \Omega$
such that $k[\Omega]\stackrel{\sim}{\to} k[X]^G$.

Consider now an object $W\in\mathcal Q_k$ with a 
$G$-invariant $\mathcal V_k$-morphism $f:X\to W$.
Then $f$ is dominant rational map
with $f^*(k[W])\subset k[X]^G$. 
Thus we obtain an $k$-algebra homomorphism $f^*:k[W]\to k[\Omega]=k[X]^G$,
which by lemma \ref{lem11} induces a $\mathcal Q_k$-morphism
from $\Omega$ to $W$.

Therefore $\Mor_{\mathcal V_k}(X,W)^G\simeq \Mor_{\mathcal Q_k}(\Omega, W)$.
\end{proof}

\section{An example}
Let $k$ be a field of characteristic zero.
In \cite{DgF} Daigle and Freudenburg gave an example of a locally nilpotent
derivation $D$ of $k[\A^5]$ such that $\ker D$ is not finitely
generated. This is the lowest-dimensional example known today.
In coordinates $x,s,t,u,v$ the derivation $D$ can be written as
\[
D = x^3 \dd s + s\dd t + t \dd u + x^2 \dd v 
\]
The associated group action of the additive group $G_a$ is given
by
\[
\mu(r):(x,s,t,u,v) \mapsto
(x,s+rx^3, 
t+rs+\frac{r^2}{2}x^3,
u+rt+\frac{r^2}{2}s+\frac{r^3}{6}x^3,
v+rx^2)
\]
The action is free outside the set of fixed points
\[
(\A^5)^\mu=\{(0,0,0,v,u):u,v\in k\}
\]
An explicit calculation shows that the following regular functions
on $\A^5$ are invariant:
\begin{alignat*}{2}
\phi_1 &= x &&\\
\phi_2 &= 2x^3t-s^2  &&\\
\phi_3 &= 3x^6u-3x^3ts + s^3 &&\\
\phi_4 &= xv-s &&\\
\phi_5 &= x^2 ts -s^2v +2x^3tv -3x^5u &\ =&\ (\phi_2\phi_4 - \phi_3)/\phi_1 \\
\phi_6 &= -18x^3tsu + 9x^6u^2 + 8x^3t^3 + 6s^3 u -3 x^6t^2s^2
       &\ = &\ (\phi_2^3 + \phi_3 ^2 ) / \phi_1^6  
\end{alignat*}

Further explicit calculations yields the following:
\begin{lemma}
Let $V=\{w=(w_1,\ldots,w_6)\in\A^6: w_5w_1=w_2w_4-w_3,
                                    w_6w_1^6=w_2^3 + w_3^2 \}$.
Then $V$ is an affine subvariety of $\A^6$ with
$\Sing V=\{w\in\A^6:w_1=w_2=w_3=0\}$.
$\Sing V$ is a Weil divisor of $V$.

The regular functions $(\phi_i)_{i=1..6}$ defined above
give an invariant morphism $\phi:\A^5\to V$ such that
\begin{enumerate}
\item
$\phi$ has rank $4$ outside $E=\{x=s=0\}$.
\item
$\phi$ maps $\A^5\setminus E$ surjectively on $V\setminus\Sing V$.
\item
For every $p\in V\setminus\Sing V$ the fiber $\phi^{-1}(p)$
coincides with a $G_a$-orbit.
\end{enumerate}
\end{lemma}
It follows that
\[
k[\A^5]^\mu = k[\A^5\setminus E]^\mu \simeq k [V\setminus\Sing V]
\]
Thus the ring of invariants is indeed isomorphic to the ring of
regular functions of a quasi affine variety, namely $V\setminus\Sing V$.

\section{Appendix}
Naturally one would prefer having a regular morphism in $(*)$ instead
of having merely a rational map. 
One might hope for a general result implying that this map is
automatically regular, and endeavor to prove a statement like the following.
\begin{Falsity}
Let $V$ , $W$ be affine varieties, $f:V\to W$ be a regular map,
$H\subset W$ a hypersurface and assume that $f^*(k[W\setminus
H])\subset k[V]$. Then $f(V)\subset W\setminus H$ (implying that one
has a regular morphism from $V$ to $W\setminus H$ and not merely a
rational map.)
\end{Falsity}
But this is wrong:
\begin{example}
Let $V=\{(x_1,x_2,x_3)\in\mathbb A^3 : x_1\ne 0\}$,
$W=\{(z_1,z_2,z_3,z_4)\in\mathbb A^4: z_1z_4=z_2z_3\}$,
$H=\{(z_1,z_2,z_3,z_4): z_1=z_2=0\}$ and $f:V\to W$ given
by $f(x_1,\ldots,x_3)=(x_1x_2,x_1x_3,x_2,x_3)$. Then $H$ is a
hypersurface in $W$, but $f^{-1}(H)=\{(x_1,0,0):x_1\ne 0\}$  is a curve, 
and therefore has
codimension two in $V$. As a consequence $f^*k[W\setminus H]
\subset k[V\setminus f^{-1}(H)]=k[V]$, although $f(V)\cap
H=\{(0,0,0,0)\}\ne\emptyset$.
\end{example}

\bigskip
\vbox{\leftskip=150pt
\parindent=0pt {\obeylines
}}
\par

\begin{thebibliography}{Bla}
\bibitem{ACN} A'Campo-Neuen, A.:
Note on a counterexample to Hilbert's
fourteenth problem.
\sl Indag. Math.(n.S.)\bf 5\rm, 253--257 (1994)
\bibitem{B} Bourbaki, N.:
Alg\'ebre Commutative. VII.
Hermann 1965.

\bibitem{DgF} Daigle, Daniel; Freudenburg, Gene:
A counterexample to Hilbert's fourteenth problem in dimension $5$. 
\sl J.~Algebra \bf 221\rm,  528--535. (1999)

\bibitem{DF} Deveney, J.K.; Finston, D.:
$G_a$-actions on $\C^3$ and $\C^7$.
\sl Comm. Alg. \bf 22\rm, 6295--6302 (1994)

\bibitem{Dx} Dixmier, J.:
Solution n\'egative du probl\`eme des invariants d'apr\'es Nagata.
\sl S\'eminaire Bourbaki, Expos\'e \bf 175\rm. Paris 1959

\bibitem{F} Freudenburg, G.:
 A counterexample to Hilbert's
   fourteenth problem in dimension six. 
\sl Transform. Groups \bf 5 \rm, no.
   1, 61--71 (2000)

\bibitem{G} Goodman, J.E.:
Affine open subsets of algebraic varieties and ample divisors.
\sl Ann. Math. \bf 89\rm, 160-183 (1969)

\bibitem{H} Hilbert, D.:
Probl\`emes future de math\'ematiques.
Compte Rendu du 2${}^{em}$ Congr\'es International des
Math\'ematiciens.
58--114.
Gauthier-Villars. Paris 1902.

\bibitem{N1} Nagata, M.: Lectures on the fourteenth problem of Hilbert.
Tata Institute, Bombay. (1965)

\bibitem{N2} Nagata, M.:
Local rings. New York (1962)

\bibitem{Ne}  Neeman, A.: Steins, affines and Hilberts fourteenth problem.
 \sl Ann. of Math. \bf 127\rm, 229--244 (1988)

\bibitem{Oe} Oeljeklaus, K.:
On the holomorphic separability of discrete quotients of complex Lie
groups.
\sl Math. Z. \bf 211\rm, (92) no. 4, 627--533

\bibitem{P} Popov, V.L.: Hilberts Theorem on Invariant.
\sl Soviet Math. Dokl. \bf 20\rm, 1318--1322 (1979)


\bibitem{R} Rees, D.: On a problem of Zariski.
 \sl Illinois J.~Math. \bf 2\rm, 145--149
(1958)

\bibitem{Rob} Roberts, P.:
An infinitely generated symbolic blow-up in a power series ring and a
new counterexample to Hilbert's fourteenth problem.
\sl J.~Alg. \bf 132\rm, 461--473 (1990)

\bibitem{Z} Zariski, O.:
Interpr\'etations algebrico-geom\'etrique du $14^{eme}$
probl\`eme de Hilbert.
\sl Bull. Sci. Math. \bf 78\rm, 155--168 (1954)


\end{thebibliography}
\end{document}